\documentclass[12pt,twoside]{article}
\usepackage{amssymb,amsmath} 
\newtheorem{theorem}{Theorem}[section]
\newtheorem{corollary}[theorem]{Corollary}
\newtheorem{lemma}[theorem]{Lemma}

\newtheorem{remark}[theorem]{Remark}

\numberwithin{equation}{section}
\newcommand \proof{\indent\textbf{Proof. \ }}
\newcommand \epf{\hfill \rule{0.2cm}{0.2cm}\bigskip \par}

\title{The solvability of Brezis-Nirenberg type problems of singular quasilinear elliptic equation}

\begin{document}
\author{{Benjin Xuan}
\thanks{Supported by Grant 10071080 and 10101024 from the NNSF of China.
}\\
{\it Department of Mathematics}\\
{\it University of Science and Technology of China}\\
{\it Universidad Nacional de Colombia}\\
{\it e-mail:wenyuanxbj@yahoo.com}}
\date{}
\maketitle

\begin{abstract}
In this paper, we consider the existence and non-existence of non-trivial solution to a Brezis-Nirenberg type problem with singular weights.
First, we obtain a compact imbedding theorem which is an extension of the classical Rellich-Kondrachov compact imbedding theorem, and consider the corresponding eigenvalue problem.
Secondly, we deduce a Pohozaev type identity and obtained a non-existence result.
Thirdly, based on a
generalized concentration compactness principle, we will give some
abstract conditions when the functional satisfies the (PS)$_c$
condition. Finally, based on the explicit form of the extremal function, we will obtain some existence results to the problem.


\end{abstract}

\section{Introduction.}\label{intro}
In this paper, we consider the existence and non-existence of non-trivial solution to the following Brezis-Nirenberg type problem with singular weights:
\begin{equation}
\label{eq1.1}
\left\{
\begin{array}{l}
-\mbox{div\,}(|x|^{-ap}|Du|^{p-2}Du)=|x|^{-bq}|u|^{q-2}u+\lambda |x|^{-(a+1)p+c}|u|^{p-2}u,\mbox{ in } \Omega\\[2mm]
u= 0, \ \ \mbox{on } \partial\Omega,
\end{array}
\right.
\end{equation}
where $\Omega\subset \mathbb{R}^n$ is an open bounded domain with $C^1$ boundary and $0\in \Omega$, $1<p<n,\ -\infty< a<\frac{n-p}p,\ a\leq b\leq a+1,\ q=p^*(a,b)=\frac{np}{n-dp},\ d=1+a-b\in [0,\ 1],\ c>0$.

The starting point of the variational approach to these problems is the following weighted Sobolev-Hardy inequality due to Caffarelli, Kohn and Nirenberg \cite{CKN}, which is called the Caffarelli-Kohn-Nirenberg inequality. Let $1<p<n$. For all $u\in C_0^\infty(\mathbb{R}^n)$, there is a constant $C_{a,b}>0$ such that
\begin{equation}
\label{eq1.2}
\Big(\int_{\mathbb{R}^n}|x|^{-bq}|u|^{q}\,dx \Big)^{p/q}\leq C_{a,b}\int_{\mathbb{R}^n}|x|^{-ap}|Du|^{p}\,dx,
\end{equation}
where
\begin{equation}
\label{eq1.3}
-\infty< a<\frac{n-p}p,\ a\leq b\leq a+1,\ q=p^*(a,b)=\frac{np}{n-dp},\ d=1+a-b.
\end{equation}

Let ${\cal D}_a^{1,p}(\Omega)$ be the completion of $C_0^\infty(\mathbb{R}^n)$, with respect to the norm $\|\cdot\|$ defined by
$$
\|u\|=\Big(\int_{\Omega}|x|^{-ap}|Du|^{p}\,dx \Big)^{1/p}.
$$
From the boundedness of $\Omega$ and the standard approximation
argument, it is easy to see that (\ref{eq1.2}) holds for any $u\in
{\cal D}_a^{1,p}(\Omega)$ in the sense:
\begin{equation}
\label{eq1.4}
\Big(\int_{\Omega}|x|^{-\alpha}|u|^{r}\,dx \Big)^{p/r}\leq C \int_{\Omega}|x|^{-ap}|Du|^{p}\,dx,
\end{equation}
for $1\leq r\leq \frac{np}{n-p},\ \frac\alpha r\leq (1+a)+n(\frac1r-\frac 1p)$, that is, the imbedding ${\cal D}_a^{1,p}(\Omega) \hookrightarrow L^r(\Omega, |x|^{-\alpha})$ is continuous, where $L^r(\Omega, |x|^{-\alpha})$ is the weighted $L^r$ space with norm:
$$
\|u\|_{r, \alpha}:=\|u\|_{L^r(\Omega, |x|^{-\alpha})}=\Big( \int_{\Omega}|x|^{-\alpha}|u|^{r}\,dx\Big)^{1/r}.
$$

On ${\cal D}_a^{1,p}(\Omega)$, we can define the energy functional:
\begin{equation}
\label{eq1.04}
E_\lambda(u)=\frac1p\int_{\Omega} |x|^{-ap}|Du|^{p}\,dx-\frac1q \int_{\Omega} |x|^{-bq}|u|^{q}\,dx-\frac\lambda p \int_{\Omega} |x|^{-(a+1)p+c}|u|^p\,dx.
\end{equation}
From (\ref{eq1.4}), $E_\lambda$ is well-defined in ${\cal D}_a^{1,p}(\Omega)$, and $E_\lambda\in C^1({\cal D}_a^{1,p}(\Omega),\mathbb{R})$. Furthermore, the critical points of $E_\lambda$ are weak solutions of problem (\ref{eq1.1}).

We note that for $p=2,\ a=b=0$ and $c=2$, problem (\ref{eq1.1}) becomes
\begin{equation}
\label{eq1.5}
\left\{
\begin{array}{l}
-\Delta u =|u|^{q-2}u+\lambda u,\mbox{ in } \Omega\\[2mm]
u= 0, \ \ \mbox{on } \partial\Omega,
\end{array}
\right.
\end{equation}
where $q=2^*=\frac{2n}{n-2}$ is the critical Sobolev exponent. Problem (\ref{eq1.5}) has been studied in a more general context in the famous paper by Brezis and Nirenberg \cite{BN}. Since the imbedding $H_0^1(\Omega) \hookrightarrow  L^q(\Omega)$ is not compact for $q=\frac{2n}{n-2}$, the corresponding energy functional does not satisfy the (PS) condition globally, which caused a serious difficulty when trying to find critical points by standard variational methods. By carefully analyzing the energy level of a cut-off function related to the extremal function of the Sobolev inequality in $\mathbb{R}^n$, Brezis and Nirenberg obtained that the energy functional does satisfy the (PS)$_c$ for some energy level $c <\frac1nS^{n/2}$, where $S$ is the best constant of the Sobolev inequality.

Brezis-Nirenberg type problems have been generalized to many situations (see \cite{CG, EH1, EH2, GV, JS, NL, PS, XC, ZXP} and references therein). In \cite{EH2, GV, ZXP}, the results of \cite{BN} had been extended to the p-Laplace case; \cite{PS, XC} extended the results of \cite{BN} to polyharmonic operators; Jannelli and Solomini \cite{JS} considered the case with singular potentials where $p=2, a=0, c=2, b\in [0,1]$; while \cite{CG} considered the weighted case where $p=2, a<\frac{n-2}2,\ b\in [a, a+1], c>0$, and \cite{NL} considered the case where $p=2,\ a=0$ and $\Omega$ is a ball.

All the above references are based on the fact that the extremal functions are symmetric and have explicit forms.
In \cite{CC}, based on a generalization of the moving plane method, Chou and Chu considered the symmetry of the extremal functions for $a\geq 0,\ p=2$; In \cite{HT}, Horiuchi successfully treated the symmetry properties of the extremal functions for the case $p>1,\ a\leq 0$ by a clever reduction to the case $a=0$ (where Schwarz symmetrization gives the symmetry of the extremal functions); On the contrary, there are some symmetry breaking results (cf. \cite{CW, BW}) for $a<0$. We define
\begin{equation}
\label{eq1.11}
S(a, b)=\inf_{u\in {\cal D}_a^{1,p}(\mathbb{R}^n)\setminus \{0\}} E_{a,b}(u),
\end{equation}
to be the best embedding constants, where
\begin{equation}
\label{eq1.12}
E_{a,b}(u)=\frac{\displaystyle \int_{\mathbb{R}^n}|x|^{-ap} |Du|^p\,dx}{\big(\displaystyle \int_{\mathbb{R}^n}|x|^{-bq} |u|^q\,dx \big)^{p/q}},
\end{equation}
and
$$
S_R(a, b)=\inf_{u\in {\cal D}_{a, R}^{1,p}(\mathbb{R}^n)\setminus \{0\}} E_{a,b}(u),
$$
where ${\cal D}_{a, R}^{1,p}(\mathbb{R}^n)=\{u\in {\cal D}_{a}^{1,p}(\mathbb{R}^n)\,|\, u \mbox{ is radial} \}$. It is well known that for $a<\frac{n-p}p$ and $b-a<1$, $S_R(a, b)$ is always achieved and the extremal functions are given by

\begin{equation}
\label{eq1.00013}
U_{a,b}(r)=c_0\Big(\frac{n-p-pa}{1+r^{\frac{dp(n-p-pa)}{(p-1)(n-dp)}}} \Big)^\frac{n-dp}{dp}
\end{equation}
where
\begin{equation}
\label{eq1.00014}
c_0=\big(\frac{n}{(p-1)^{p-1}(n-dp)}\big)^\frac{n-dp}{dp^2}.
\end{equation}
Under some condition on parameters $a,\ b,\ n,\ p$, \cite{CW, BW} obtain that $S(a,b)<S_R(a,b)$ for $a<0$. In this case, it is very difficult to verify that the corresponding energy functional satisfies the (PS)$_c$ condition.

In section 2, based on the Caffarelli-Kohn-Nirenberg inequality and the classical Rellich-Kondrachov compactness theorem, we will first deduce a compact imbedding theorem and study the corresponding eigenvalue problem:
\begin{equation}
\label{eq1.6}
\left\{
\begin{array}{l}
-\mbox{div\,}(|x|^{-ap}|Du|^{p-2}Du)=\lambda |x|^{-(a+1)p+c}|u|^{p-2}u,\mbox{ in } \Omega\\[2mm]
u= 0, \ \ \mbox{on } \partial\Omega.
\end{array}
\right.
\end{equation}
In section 3, based on a Pohozaev type
identity, we obtained a non-existence result for problem
(\ref{eq1.1}) with $\lambda\leq 0$. In section 4, based on a
generalized concentration compactness principle, we shall give some abstract conditions when the functional satisfies the (PS)$_c$
condition. In section 5, based on the explicit form of the extremal function, we will obtain some existence results
to problem (\ref{eq1.1}).

\section{Eigenvalue problem in domain general}
In this section, we first deduce a compact imbedding theorem which
is an extension of the classical Rellich-Kondrachov compactness
theorem.
\begin{theorem}[Compact imbedding theorem]
\label{thm2.1} Suppose that $\Omega\subset \mathbb{R}^n$ is an
open bounded domain with $C^1$ boundary and $0\in \Omega$,
$1<p<n,\ -\infty< a<\frac{n-p}p$. The imbedding ${\cal
D}_a^{1,p}(\Omega) \hookrightarrow L^r(\Omega, |x|^{-\alpha})$ is
compact if $1\leq r< \frac{np}{n-p},\ \alpha < (1+a)r+n(1-\frac
rp)$.
\end{theorem}
\proof The continuity of the imbedding is a direct consequence of
the Caffarelli-Kohn-Nirenberg inequality (\ref{eq1.2}) or
(\ref{eq1.4}). To prove the compactness, let $\{u_m\}$ be a
bounded sequence in ${\cal D}_a^{1,p}(\Omega)$. For any $\rho>0$
with $B_\rho(0)\subset \Omega$ is a ball centered at the origin
with radius $\rho$, there holds $\{u_m\}\subset
W^{1,p}(\Omega\setminus B_\rho(0))$. Then the classical
Rellich-Kondrachov compactness theorem guarantees the existence of
a convergent subsequence of $\{u_m\}$ in $L^r(\Omega\setminus
B_\rho(0))$. By taking a diagonal sequence, we can assume without
loss of generality that $\{u_m\}$ converges in
$L^r(\Omega\setminus B_\rho(0))$ for any $\rho>0$.

On the other hand, for any $1\leq r< \frac{np}{n-p}$, there exists a $b\in (a, a+1]$ such that $r<q=p^*(a, b)=\frac{np}{n-dp},\ d=1+a-b\in [0,\ 1)$. From the Caffarelli-Kohn-Nirenberg inequality (\ref{eq1.2}) or (\ref{eq1.4}), $\{u_m\}$ is also bounded in $L^q(\Omega, |x|^{-bq})$. By the H\"{o}der inequality, for any $\delta>0$, there holds
\begin{equation}
\label{eq2.1}
\begin{array}{ll}
&\displaystyle \int_{|x|<\delta}|x|^{-\alpha}|u_m-u_j|^{r}\,dx \\[3mm]
&\ \ \ \ \  \leq \Big( \displaystyle \int_{|x|<\delta}|x|^{-(\alpha-br)\frac q{q-r}}\,dx\Big)^{1-\frac rq}
 \Big(\displaystyle \int_{\Omega}|x|^{-br}|u_m-u_j|^{r}\,dx\Big)^{r/q}\\[3mm]
&\ \ \ \ \   \leq C \Big( \displaystyle \int_0^\delta r^{n-1-(\alpha-br)\frac q{q-r}}\,dr\Big)^{1-\frac rq}\\[3mm]
&\ \ \ \ \  =C \delta ^{n-(\alpha-br)\frac q{q-r}},
\end{array}
\end{equation}
where $C>0$ is a constant independent of $m$. Since $\alpha
< (1+a)r+n(1-\frac rp)$, there holds $n-(\alpha-br)\frac
q{q-r}>0$. Therefore, for a given $\varepsilon>0$, we first fix
$\delta>0$ such that
$$
\int_{|x|<\delta}|x|^{-\alpha}|u_m-u_j|^{r}\,dx \leq \frac{\varepsilon}2, \ \forall\ m, j\in \mathbb{N}.
$$
Then we choose $N\in  \mathbb{N}$ such that
$$
\int_{\Omega\setminus B_\delta(0)}|x|^{-\alpha}|u_m-u_j|^{r}\,dx \leq C_\alpha \int_{\Omega\setminus B_\delta(0)} |u_m-u_j|^{r}\,dx \leq \frac{\varepsilon}2, \ \forall\ m, j\geq N,
$$
where $C_\alpha=\delta ^{-\alpha}$ if $\alpha\geq 0$ and $C_\alpha=(\mbox{diam\,}(\Omega) )^{-\alpha}$ if $\alpha< 0$. Thus
$$
\int_{\Omega}|x|^{-\alpha}|u_m-u_j|^{r}\,dx \leq \varepsilon , \ \forall\ m, j\geq N,
$$
that is, $\{u_m\}$ is a Cauchy sequence in $L^q(\Omega, |x|^{-bq})$.
\epf

\begin{remark}
\label{rem2.2}
\cite{CC} had obtained Theorem \ref{thm2.1} for the case $p=2$.
\end{remark}

In order to study the eigenvalue problem (\ref{eq1.6}), let us introduce the following functionals in ${\cal D}_a^{1,p}(\Omega)$:
$$
\Phi(u):= \int_{\Omega} |x|^{-ap}|Du|^{p}\,dx,\ \mbox{and } J(u):=\int_{\Omega} |x|^{-(a+1)p+c}|u|^p\,dx.
$$
For $c>0$, $J$ is well-defined. Furthermore, $\Phi, J\in C^1({\cal D}_a^{1,p}(\Omega),\mathbb{R})$, and a real value $\lambda$ is an eigenvalue of problem (\ref{eq1.6}) if and only if there exists $u\in {\cal D}_a^{1,p}(\Omega)\setminus\{0\}$ such that $\Phi^\prime(u)=\lambda J^\prime(u)$. At this point let us introduce set
$$
{\cal M}:=\{u\in {\cal D}_a^{1,p}(\Omega)\ :\ J(u)=1 \}.
$$
Then ${\cal M}\neq \emptyset$ and ${\cal M}$ is a $C^1$ manifold in ${\cal D}_a^{1,p}(\Omega)$. It follows from the standard  arguments that eigenvalues of (\ref{eq1.6}) correspond to critical values of $\Phi|_{{\cal M}}$. From Theorem \ref{thm2.1}, $\Phi$ satisfies the (PS) condition on ${\cal M}$. Thus a sequence of critical values of $\Phi|_{{\cal M}}$ comes from the Ljusternik-Schnirelman critical point theory on $C^1$ manifolds. Let $\gamma(A)$ denote the Krasnoselski's genus on ${\cal D}_a^{1,p}(\Omega)$ and for any $k\in \mathbb{N}$, set
$$
\Gamma_k:=\{A\subset {\cal M}\ :\ A \mbox{ is compact, symmetric and } \gamma(A)\geq k\}.
$$
Then values
\begin{equation}
\label{eq2.2}
\lambda_k:=\inf_{A\in\Gamma_k }\max_{u\in A} \Phi(u)
\end{equation}
are critical values and thence are eigenvalues of problem (\ref{eq1.6}). Moreover, $\lambda_1\leq \lambda_2\leq \cdots \leq   \lambda_k\leq \cdots \to+\infty$.

From the Caffarelli-Kohn-Nirenberg inequality (\ref{eq1.2}) or (\ref{eq1.4}), it is easy to see that
$$
\lambda_1=\inf\{\Phi(u)\ :\ u\in {\cal D}_a^{1,p}(\Omega),  J(u)=1\}>0,
$$
and the corresponding eigenfunction $e_1\geq 0$.

\section{Pohozaev identity and non-existence result}
In this section, we deduce a Pohozaev-type identity and obtain some non-existence results. First let us recall the following Pohozaev integral identity due to Pucci and Serrin \cite{PS1}:
\begin{lemma}
\label{lem3.1}
Let $u\in C^2(\Omega)\cap C^1(\bar \Omega)$ be a solution of the Euler-Lagrange equation
\begin{equation}
\label{eq3.1}
\left\{
\begin{array}{l}
\mbox{div\,}\{{\cal F}_p(x,u,Du)\}={\cal F}_u(x,u, Du),\mbox{ in } \Omega\\[2mm]
u= 0, \ \ \mbox{on } \partial\Omega,
\end{array}
\right.
\end{equation}
where $p=(p_1,\cdots, p_n)=Du=(\partial u/\partial x_1,\cdots, \partial u/\partial x_n)$ and ${\cal F}_u=\partial {\cal F}/\partial u$. Let $A$ and $h$ be, respectively, scalar and vector-value function of class $C^1(\Omega)\cap C(\bar \Omega)$. Then there holds
\begin{equation}
\label{eq3.2}
\begin{array}{ll}
&\displaystyle \oint_{\partial \Omega}\Big[{\cal F}(x,0,Du)-\dfrac{\partial u}{\partial x_i}{\cal F}_{p_i}(x,0,Du)  \Big](h\cdot \nu)\,ds \\[3mm]
& \ \ \ \ \ \ =\displaystyle \int_{\Omega}\Big\{ {\cal F}(x,u,Du)\mbox{div\,}h+ h_i{\cal F}_{x_i}(x,u,Du)\\[3mm]
& \ \ \ \ \  \ \ \ \ \ \ \ \ -\big[ \dfrac{\partial u}{\partial x_j} \dfrac{\partial h_j}{\partial x_i} + u\dfrac{\partial A}{\partial x_i}\big]{\cal F}_{p_i}(x,u,Du)\\[3mm]
& \ \ \ \ \  \ \ \ \ \ \  \ \ -A \big[\dfrac{\partial u}{\partial x_i}{\cal F}_{p_i}(x,u,Du)+u {\cal F}_{u}(x,u,Du) \big]
\Big\}\,dx,
\end{array}
\end{equation}
where repeated indices $i$ and $j$ are understood to be summed from $1$ to $n$.
\end{lemma}

Let us consider the following problem:
\begin{equation}
\label{eq3.3}
\left\{
\begin{array}{l}
-\mbox{div\,}(|x|^{-ap}|Du|^{p-2}Du)=g(x,u),\mbox{ in } \Omega\\[2mm]
u= 0, \ \ \mbox{on } \partial\Omega,
\end{array}
\right.
\end{equation}
where $g$ satisfies $g(x, 0)=0$.
Suppose that ${\cal F}(x,u, Du)=\frac{1}{p}|x|^{-ap}|Du|^p-G(x,u)$, where $G(x,u)=\int_0^u g(x,t)\,dt$ is the primitive of $g(x,u)$. If we choose $h(x)=x,\ A=\frac np-(1+a)$, then (\ref{eq3.2}) becomes
\begin{equation}
\label{eq3.4}
\begin{array}{ll}
&(1-\dfrac1p)\displaystyle \oint_{\partial \Omega} (x\cdot \nu)|\dfrac{\partial u}{\partial \nu}|^p\,ds \\[3mm]
& \ \ \ \ \ \ =\displaystyle \int_{\Omega}\big[nG(x,u)+(x, G_x)+(1+a-\dfrac np)ug(x,u) \big]\,dx.
\end{array}
\end{equation}
As to problem (\ref{eq1.1}), suppose that $G(x,u)=\dfrac1q|x|^{-bq}|u|^q+\dfrac\lambda p |x|^{-p(1+a)+c}|u|^p$, then (\ref{eq3.2}) or (\ref{eq3.4}) becomes
\begin{equation}
\label{eq3.5}
(1-\frac1p)\displaystyle \oint_{\partial \Omega} (x\cdot \nu)|\dfrac{\partial u}{\partial \nu}|^p\,ds =\frac{c\lambda}p \int_{\Omega} |x|^{-(a+1)p+c}|u|^p\,dx.
\end{equation}

Thus we obtain the following non-existence result:
\begin{theorem}
\label{thm3.2}
There is no solution to problem (\ref{eq1.1}) when $\lambda\leq 0$ and $\Omega$ is a (smooth) star-shaped domain with respect to the origin.
\end{theorem}
\proof The above deduction is formal. In fact, the solution to
problem (\ref{eq1.1}) may not be of class $C^2(\Omega)\cap
C^1(\bar \Omega)$. We need the approximation arguments
in \cite{GV} and \cite{CG} (cf. Appendix A). \epf

\section{(PS)$_c$ condition}
In this section, we first give a concentration compactness principle which is a weighted version of the Concentration Compactness Principle II due to P. L. Lions \cite{LPL1, LPL2}.

\begin{theorem}[Concentration Compactness Principle]
\label{thm4.1}
Let $1<p<n,\ -\infty< a<\frac{n-p}p,\ a\leq b\leq a+1,\ q=p^*(a,b)=\frac{np}{n-dp},\ d=1+a-b\in [0,\ 1]$, and ${\cal M}(\mathbb{R}^n)$ be the space of bounded measures on $\mathbb{R}^n$. Suppose that $\{u_m\}\subset {\cal D}_a^{1,p}(\mathbb{R}^n)$ be a sequence such that:
$$
\begin{array}{ll}
u_m\rightharpoonup u & \mbox{ in }   {\cal D}_{a}^{1,p}(\mathbb{R}^n),\\[2mm]
\mu_m:=\big||x|^a Du_m| \big|^p\,dx \rightharpoonup \mu & \mbox{ in }   {\cal M}(\mathbb{R}^n),\\[2mm]
\nu_m:=\big||x|^b u_m| \big|^q\,dx  \rightharpoonup \nu & \mbox{ in }   {\cal M}(\mathbb{R}^n),\\[2mm]
u_m\to u & \mbox{ a.e. on }   \mathbb{R}^n.
\end{array}
$$
There there hold:
\begin{enumerate}
\item[(1)] There exists some at most countable set $J$, a family $\{x^{(j)}\ :\ j\in J\}$ of distinct points in $\mathbb{R}^n$, and a family $\{\nu^{(j)}\ :\ j\in J \}$ of positive numbers such that
\begin{equation}
\label{eq4.1}
\nu=\big||x|^{-b} u| \big|^q\,dx+\sum_{j\in J} \nu^{(j)}\delta_{x^{(j)}},
\end{equation}
where $\delta_x$ is the Dirac-mass of mass $1$ concentrated at $x\in \mathbb{R}^n$.

\item[(2)] There holds
\begin{equation}
\label{eq4.2}
\mu\geq \big||x|^{-a} Du| \big|^p\,dx+\sum_{j\in J} \mu^{(j)}\delta_{x^{(j)}},
\end{equation}
for some family $\{\mu^{(j)}>0\ :\ j\in J \}$ satisfying
\begin{equation}
\label{eq4.3}
S(a,b)\big( \nu^{(j)}\big)^{p/q}\leq \mu^{(j)},\ \ \mbox{for all }j\in J.
\end{equation}
In particular, $\sum\limits_{j\in J}\big( \nu^{(j)}\big)^{p/q}<\infty$.
\end{enumerate}
\end{theorem}
\proof The proof is similar to that of the Concentration
Compactness Principle II (see also \cite{SM}).

\begin{theorem}
\label{thm4.2}
Let $1<p<n,\ -\infty< a<\frac{n-p}p,\ a\leq b\leq a+1,\ q=p^*(a,b)=\frac{np}{n-dp},\ d=1+a-b\in (0,\ 1],\ c>0$ and $0<\lambda<\lambda_1$. Then functional $E_\lambda$ defined in (\ref{eq1.04}) satisfies the (PS)$_c$ condition in ${\cal D}_{a}^{1,p}(\Omega)$ at the energy level $M< \frac dn S(a,b)^{\frac n{dp}}$.
\end{theorem}
\indent\textbf{Proof. \ 1. }The boundedness of (PS)$_c$ sequence.

Suppose that $\{u_m\}\subset {\cal D}_{a}^{1,p}(\Omega)$ is a (PS)$_c$ sequence of functional $E_\lambda$, that is,
$$
E_\lambda(u_m)\to M,\  \mbox{and } E_\lambda^\prime(u_m)\to 0, \ \mbox{in } ({\cal D}_{a}^{1,p}(\Omega))^\prime.
$$
Then as $m\to \infty$, there hold
\begin{equation}
\label{eq4.7}
\begin{array}{ll}
&M+o(1)=E_\lambda(u_m)\\[3mm]
&\ \ \  =\dfrac1p\displaystyle \int_{\Omega} |x|^{-ap}|Du_m|^{p}\,dx-\dfrac1q \displaystyle \int_{\Omega} |x|^{-bq}|u_m|^{q}\,dx-\dfrac\lambda p \displaystyle \int_{\Omega} |x|^{-(a+1)p+c}|u_m|^p\,dx
\end{array}
\end{equation}
and
\begin{equation}
\label{eq4.8}
\begin{array}{ll}
 o(1)\|\varphi\| & =(E_\lambda(u_m), \varphi) \\[3mm]
&=\displaystyle \int_{\Omega} |x|^{-ap}|Du_m|^{p-2}Du_m\cdot D\varphi\,dx- \displaystyle \int_{\Omega} |x|^{-bq}|u_m|^{q-2}u_m\varphi\,dx\\[3mm]
&\ \ \ \ \ \ \ \ -\displaystyle \lambda \int_{\Omega} |x|^{-(a+1)p+c}|u_m|^{p-2}u_m\varphi\,dx
\end{array}
\end{equation}
for any $\varphi\in {\cal D}_{a}^{1,p}(\Omega)$, where $o(1)$ denotes any quantity that tends to zero as $m\to \infty$.
From (\ref{eq4.7}) and (\ref{eq4.8}), $m\to \infty$, there holds
\begin{equation}
\label{eq4.9}
\begin{array}{ll}
 qM+o(1)+o(1)\|u_m\|& =qE_\lambda(u_m)-(E_\lambda(u_m), v)\\[2mm]
&=(\dfrac qp-1)\displaystyle \int_{\Omega} |x|^{-ap}|Du_m|^{p}\,dx\\[3mm]
&\ \ \ \ -\lambda (\dfrac qp-1)\int_{\Omega} |x|^{-(a+1)p+c}|u_m|^{p-2}u_mv\,dx\\[3mm]
&=(\dfrac qp-1)(1-\dfrac{\lambda}{\lambda_1})\|u_m\|^p,
\end{array}
\end{equation}
that is, $\{u_m\}$ is bounded in ${\cal D}_{a}^{1,p}(\Omega)$, since $q>p, \lambda<\lambda_1$. Thus up to a subsequence, there hold
$$
\begin{array}{ll}
u_m\rightharpoonup u & \mbox{ in }   {\cal D}_{a}^{1,p}(\Omega),\\[2mm]
u_m\rightharpoonup u & \mbox{ in }   L^{q}(\Omega, |x|^{-bq}),\\[2mm]
u_m\to u & \mbox{ in }   L^{r}(\Omega, |x|^{-\alpha}), \ \forall\  1\leq r< \frac{np}{n-p},\ \frac\alpha r< (1+a)+n(\frac1r-\frac 1p) \\[2mm]
u_m\to u & \mbox{ a.e. on }   \Omega.
\end{array}
$$
From the concentration compactness principle-Theorem \ref{thm4.1}, there exist non-negative measures $\mu, \nu$ and a countable family $\{x_j\}\subset \bar \Omega$ such that
$$
\begin{array}{ll}
&|x|^{-b}|u_m|^{q}\,dx \rightharpoonup \nu=\big||x|^{-b} u| \big|^q\,dx+\sum\limits_{j\in J} \nu^{(j)}\delta_{x^{(j)}},\\[2mm]
&\big||x|^{-a} Du_m| \big|^p\,dx\rightharpoonup \mu\geq \big||x|^{-a} Du| \big|^p\,dx+S(a,b)\sum\limits_{j\in J}\big( \nu^{(j)}\big)^{p/q} \delta_{x^{(j)}}.
\end{array}
$$

\indent\textbf{2. }Up to a subsequence, $u_m\to u$ in $L^q(\Omega, |x|^{-bq})$.

Since $\{u_m\}$ is bounded in ${\cal D}_{a}^{1,p}(\Omega)$, we may suppose, without loss of generality, that there exists $T\in \Big(L^{p^\prime}(\Omega, |x|^{-ap})\Big)^{n}$ such that
$$
|Du_m|^{p-2}Du_m\rightharpoonup T \mbox{ in }  \Big(L^{p^\prime}(\Omega, |x|^{-ap})\Big)^{n}.
$$

On the other hand, $|u_m|^{q-2}u_m$ is also bounded in $L^{q^\prime}(\Omega, |x|^{-bq})$ and
$$
|u_m|^{q-2}u_m\rightharpoonup |u|^{q-2}u \mbox{ in }  L^{q^\prime}(\Omega, |x|^{-bq})
$$
Taking $m\to \infty$ in (\ref{eq4.8}), there holds
\begin{equation}
\label{eq4.10}
\int_{\Omega} |x|^{-ap}T \cdot D\varphi\,dx= \int_{\Omega} |x|^{-bq}|u|^{q-2}u\varphi\,dx
+\lambda \int_{\Omega} |x|^{-(a+1)p+c}|u|^{p-2}u\varphi\,dx,
\end{equation}
for any $\varphi\in {\cal D}_{a}^{1,p}(\Omega)$.

Let $\varphi=\psi u_m$ in (\ref{eq4.8}), where $\psi\in C(\bar \Omega)$, then there holds
\begin{equation}
\label{eq4.11}
\begin{array}{ll}
& \displaystyle \int_{\Omega} |x|^{-ap}|Du_m|^{p-2}Du_m\cdot D\varphi\,dx= \displaystyle \int_{\Omega} |x|^{-bq}|u_m|^{q-2}u_m\varphi\,dx\\[3mm]
&\ \ \ \ \ \ \ \ \ +\lambda \displaystyle \int_{\Omega} |x|^{-(a+1)p+c}|u_m|^{p-2}u_m\varphi+o(1).
\end{array}
\end{equation}
Taking $m\to \infty$ in (\ref{eq4.11}), there holds
\begin{equation}
\label{eq4.12}
\int_{\Omega} \psi\, d\mu+ \int_{\Omega} |x|^{-ap}uT \cdot D\psi\,dx= \int_{\Omega} \psi\,d\nu+\lambda \int_{\Omega} |x|^{-(a+1)p+c}|u|^{p}\psi\,dx.
\end{equation}
Let $\varphi=\psi u$ in (\ref{eq4.10}), then there holds
\begin{equation}
\label{eq4.13}
\begin{array}{ll}
& \displaystyle\int_{\Omega} |x|^{-ap}uT \cdot D\psi\,dx+\displaystyle\int_{\Omega} |x|^{-ap}\psi T \cdot D u \,dx \\[3mm]
& \ \ \ \ \ \ = \displaystyle\int_{\Omega} |x|^{-bq}|u|^{q}\psi\,dx+\lambda\displaystyle \int_{\Omega} |x|^{-(a+1)p+c}|u|^{p}\psi\,dx.
\end{array}
\end{equation}
Thus (\ref{eq4.12})$-$(\ref{eq4.13}) implies that
\begin{equation}
\label{eq4.14}
\int_{\Omega} \psi\, d\mu= \sum_{j\in J}\nu_j\psi(x_j)+\int_{\Omega} |x|^{-ap}\psi T \cdot D u \,dx,
\end{equation}
which implies that
$$
S(a,b)\big( \nu^{(j)}\big)^{p/q}\leq \mu(x_j)=\nu_j.
$$
Thence $\nu_j\geq S(a,b)^{\frac n{dp}}$ if $\nu_j\neq 0$.

On the other hand, from (\ref{eq4.7}), (\ref{eq4.10}) and (\ref{eq4.14}), there holds
\begin{equation}
\label{eq4.15}
\begin{array}{ll}
M& = \dfrac1p\displaystyle\int_{\Omega} \,d\mu -\dfrac1q\displaystyle\int_{\Omega}\,d\nu -\dfrac{\lambda}p \int_{\Omega} |x|^{-(a+1)p+c}|u|^{p}\,dx \\[3mm]
& =\dfrac1p  \displaystyle\sum\limits_{j\in J}\nu_j + \dfrac1p \displaystyle\int_{\Omega} |x|^{-ap}T\cdot Du\,dx-\dfrac1q \sum\limits_{j\in J}\nu_j -\dfrac1q\displaystyle\int_{\Omega} |x|^{-bq}|u|^{q}\,dx \\[3mm]
& \ \ \  \ \ \  -\dfrac{\lambda}p \displaystyle\int_{\Omega} |x|^{-(a+1)p+c}|u|^{p}\,dx \\[3mm]
& = (\dfrac1p -\dfrac1q)\displaystyle\sum\limits_{j\in J}\nu_j+(\dfrac1p -\dfrac1q)\displaystyle\int_{\Omega} |x|^{-bq}|u|^{q}\,dx  \\[3mm]
&\geq (\dfrac1p -\dfrac1q)\displaystyle\sum\limits_{j\in J}\nu_j=\dfrac dn \displaystyle\sum\limits_{j\in J}\nu_j.
\end{array}
\end{equation}
Since it has been shown that $\nu_j\geq S(a,b)^{\frac n{dp}}$ if $\nu_j\neq 0$, the condition $c< \frac dn S(a,b)^{\frac n{dp}}$ implies that $\nu_j= 0$ for all $j\in J$. Hence there holds
$$
\int_{\Omega} |x|^{-bq}|u_m|^{q}\,dx\to \int_{\Omega} |x|^{-bq}|u|^{q}\,dx.
$$
Thus the Brezis-Lieb Lemma \cite{BL} implies that $u_m\to u$ in $L^q(\Omega, |x|^{-bq})$.

\indent\textbf{3. }Existence of convergent subsequence.

To show that $u_m\to u$ in ${\cal D}_{a}^{1,p}(\Omega)$, from the Brezis-Lieb Lemma \cite{BL}, it suffices to show that $Du_m \to Du$ a.e. in $\Omega$ and $\|u_m\|\to \|u\|$.

To show that $Du_m \to Du$ a.e. in $\Omega$, first note that
\begin{equation}
\label{eq4.18}
|x|^{-ap}( |Du_m|^{p-2}Du_m - |Du|^{p-2}Du)\cdot (Du_m-Du)\geq 0,
\end{equation}
the equality holds if and only if $Du_m=Du$.

Secondly, let $\varphi=u_m$ and $\varphi=u$ in (\ref{eq4.8}) and then let $m\to \infty$, respectively, there hold
\begin{equation}
\label{eq4.16}
\begin{array}{ll}
\|u_m\|^p&=\displaystyle \int_{\Omega} |x|^{-ap}|Du_m|^{p}\,dx\\[3mm]
&=\displaystyle \int_{\Omega} |x|^{-bq}|u_m|^{q}\,dx-\lambda\displaystyle \int_{\Omega} |x|^{-(a+1)p+c}|u_m|^{p}\,dx+o(1)\|u_m\|\\[3mm]
& \to  \displaystyle\int_{\Omega} |x|^{-bq}|u|^{q}\,dx-\lambda\displaystyle \int_{\Omega} |x|^{-(a+1)p+c}|u|^{p}\,dx
\end{array}
\end{equation}
and
\begin{equation}
\label{eq4.17}
\begin{array}{ll}
&\displaystyle \int_{\Omega} |x|^{-ap}|Du_m|^{p-2}Du_m\cdot Du \,dx\\[3mm]
& \ \ \ \  = \displaystyle\int_{\Omega} |x|^{-bq}|u_m|^{q-2}u_m u\,dx-\lambda \displaystyle\int_{\Omega} |x|^{-(a+1)p+c}|u_m|^{p-2}u_m u\,dx +o(1)\|u\|\\[3mm]
& \ \ \ \   \to \displaystyle \int_{\Omega} |x|^{-bq}|u|^{q}\,dx-\lambda \displaystyle\int_{\Omega} |x|^{-(a+1)p+c}|u|^{p}\,dx.
\end{array}
\end{equation}
From (\ref{eq4.16}) and (\ref{eq4.17}), there holds
\begin{equation}
\label{eq4.19}
\begin{array}{ll}
& \displaystyle \int_{\Omega} |x|^{-ap}( |Du_m|^{p-2}Du_m - |Du|^{p-2}Du)\cdot (Du_m-Du) \,dx\\[3mm]
& \ \ \ \ \ =\displaystyle \int_{\Omega} |x|^{-ap}|Du_m|^{p}\,dx - \displaystyle \int_{\Omega} |x|^{-ap} |Du_m|^{p-2}Du_m \cdot Du \,dx\\[3mm]
& \ \ \ \ \  \ \ \ \ \ \
-\displaystyle \int_{\Omega} |x|^{-ap} |Du|^{p-2}Du\cdot (Du_m-Du) \,dx\\[3mm]
&\ \ \ \ \  \to  0.
\end{array}
\end{equation}
(\ref{eq4.18}) and (\ref{eq4.19}) imply that $Du_m \to Du$ a.e. in $\Omega$, hence $T=|Du|^{p-2}Du$, that is, $|Du_m|^{p-2}Du_m \rightharpoonup |Du|^{p-2}Du$ in $\Big(L^{p^\prime}(\Omega, |x|^{-ap})\Big)^{n}$.

To show that $\|u_m\|\to \|u\|$, from (\ref{eq4.16}) and (\ref{eq4.17}), there holds
$$
\begin{array}{ll}
\|u\|^p & \gets \displaystyle \int_{\Omega} |x|^{-ap}|Du_m|^{p-2}Du_m\cdot Du \,dx\\[3mm]
& \ \ \ \  = \displaystyle\int_{\Omega} |x|^{-bq}|u_m|^{q-2}u_m u\,dx-\lambda \displaystyle\int_{\Omega} |x|^{-(a+1)p+c}|u_m|^{p-2}u_m u\,dx\\[3mm]
& \ \ \ \   \to \displaystyle \int_{\Omega} |x|^{-bq}|u|^{q}\,dx-\lambda \displaystyle\int_{\Omega} |x|^{-(a+1)p+c}|u|^{p}\,dx,
\end{array}
$$
thus, $\|u_m\|^p\to \|u\|^p$.
\epf

As indicated in the introduction, for $a<0$, $S(a,b)<S_R(a,b)$ and
there is no explicit form of the minimizers of $S(a,b)$, so it is
difficult to show that there exists a minimax value $M< \frac dn S(a,b)^{\frac n{dp}}$. But there does exist an
explicit form of the extremal functions of $S_R(a,b)$, the method in
\cite{BN} can be used to show that there exists a minimax value
$M< \frac dn S_R(a,b)^{\frac n{dp}}$. Next theorem
shows that in the space of radial functions, the functional
$E_\lambda$ defined in (\ref{eq1.04}) satisfies the (PS)$_c$
condition in ${\cal D}_{a,R}^{1,p}(\Omega)$ at the energy level
$M< \frac dn S_R(a,b)^{\frac n{dp}}$ in the case $p=2$.
\begin{theorem}
\label{thm4.3}
Let $\Omega=B_1(0)$ the unit ball in $\mathbb{R}^n$, $p=2<n,\ -\infty< a<\frac{n-2}2,\ a\leq b\leq a+1,\ q=2^*(a,b)=\frac{2n}{n-2d},\ d=1+a-b\in [0,\ 1],\ c>0$ and $0<\lambda<\lambda_{1}$. Then functional $E_\lambda$ defined in (\ref{eq1.04}) satisfies the (PS)$_c$ condition in ${\cal D}_{a,R}^{1,2}(\Omega)$ at the energy level $M< \frac dn S_R(a,b)^{\frac n{2d}}$.
\end{theorem}
\indent\textbf{Proof. 1. }As in the proof of Theorem \ref{thm4.2}, any (PS)$_c$ sequence is bounded in ${\cal D}_{a,R}^{1,2}(\Omega)$, and up to a subsequence, there hold
$$
\begin{array}{ll}
u_m\rightharpoonup u & \mbox{ in }   {\cal D}_{a,R}^{1,2}(\Omega),\\[2mm]
u_m\rightharpoonup u & \mbox{ in }   L^{q}(\Omega, |x|^{-bq}),\\[2mm]
u_m\to u & \mbox{ in }   L^{r}(\Omega, |x|^{-\alpha}), \ \forall\  1\leq r< \frac{2n}{n-2},\ \frac\alpha r< (1+a)+n(\frac1r-\frac 12) \\[2mm]
u_m\to u & \mbox{ a.e. on }   \Omega.
\end{array}
$$
Thence $u$ satisfies the following equation in weak sense
\begin{equation}
\label{eq4.20}
\left\{
\begin{array}{l}
-\mbox{div\,}(|x|^{-2a} Du)=|x|^{-bq}|u|^{q-2}u+\lambda |x|^{-2(a+1)+c}u,\mbox{ in } \Omega\\[2mm]
u= 0, \ \ \mbox{on } \partial\Omega.
\end{array}
\right.
\end{equation}
Thus there holds
\begin{equation}
\label{eq4.21}
\begin{array}{ll}
E_\lambda(u)&=\dfrac12 \displaystyle \int_{\Omega} |x|^{-2a}|Du|^2\,dx -\dfrac1q \displaystyle \int_{\Omega} |x|^{-bq} |u|^{q} \,dx-\dfrac \lambda 2 \displaystyle \int_{\Omega} |x|^{-2(a+1)+c}u^2\,dx \\[3mm]
&=(\dfrac12-\dfrac1q)( \displaystyle \int_{\Omega} |x|^{-2a}|Du|^2\,dx -\lambda \displaystyle \int_{\Omega} |x|^{-2(a+1)+c}u^2\,dx )  \geq0.
\end{array}
\end{equation}

\indent\textbf{2. } Let $v_m:=u_m-u$, the Brezis-Lieb Lemma \cite{BL} leads to
$$
\int_{\Omega} |x|^{-bq} |u_m|^{q} \,dx=\int_{\Omega} |x|^{-bq} |u|^{q} \,dx+\int_{\Omega} |x|^{-bq} |v_m|^{q} \,dx+o(1).
$$
From $E_\lambda(u_m)\to c$ and $(E_\lambda^\prime(u_m),u_m)\to 0$, there hold
\begin{equation}
\label{eq4.22}
\begin{array}{ll}
E_\lambda(u_m)&=E_\lambda(u) + \dfrac12 \displaystyle \int_{\Omega} |x|^{-2a}|Dv_m|^2\,dx \\[3mm]
& \ \ \ \ \ \ \ \ \  -\dfrac1q \displaystyle \int_{\Omega} |x|^{-bq} |v_m|^{q} \,dx-\dfrac \lambda 2 \displaystyle \int_{\Omega} |x|^{-2(a+1)+c}v_m^2\,dx \\[3mm]
&\to M
\end{array}
\end{equation}
and
\begin{equation}
\label{eq4.23}
\begin{array}{ll}
& \displaystyle \int_{\Omega} |x|^{-2a}|Dv_m|^2\,dx
-\displaystyle \int_{\Omega} |x|^{-bq} |v_m|^{q} \,dx- \lambda \displaystyle \int_{\Omega} |x|^{-2(a+1)+c}v_m^2\,dx \\[3mm]
& \ \ \ \ \ \to \displaystyle \int_{\Omega} |x|^{-bq} |u|^{q} \,dx+\lambda \displaystyle \int_{\Omega} |x|^{-2(a+1)+c}u^2\,dx- \displaystyle \int_{\Omega} |x|^{-2a}|Du|^2\,dx \\[3mm]
& \ \ \ \ \ =-(E_\lambda^\prime(u),u)=0.
\end{array}
\end{equation}

Up to a subsequence, we may assume that
$$
\int_{\Omega} |x|^{-2a}|Dv_m|^2\,dx- \lambda \int_{\Omega} |x|^{-2(a+1)+c}v_m^2\,dx \to b, \ \int_{\Omega} |x|^{-bq} |v_m|^{q} \,dx\to b,
$$
for some $b\geq 0$. From Theorem \ref{thm2.1}, $v_m\to 0$ in $L^2(\Omega, |x|^{-2(a+1)+c})$, then
$$
\int_{\Omega} |x|^{-2a}|Dv_m|^2\,dx\to b.
$$
On the other hand, there holds
$$
\int_{\Omega} |x|^{-2a}|Dv_m|^2\,dx \geq S_R(a,b) \Big( \int_{\Omega} |x|^{-bq} |v_m|^{q} \,dx\Big)^{2/q}.
$$
Thus there holds $b\geq S_R(a,b) b^{2/q}$, either $b\geq S_R(a,b)^{\frac n{2d}}$ or $b=0$. If $b=0$, the proof is complete. Assume that $b\geq S_R(a,b)^{\frac n{2d}}$, from (\ref{eq4.21}) and (\ref{eq4.22}), there holds
$$
\frac dn S_R(a,b)^{\frac n{2d}}\leq (\frac12-\frac1q)b\leq M < \frac dn S_R(a,b)^{\frac n{2d}}
$$
a contradiction. \epf

\section{Existence results}
In this section, by verifying that there exists a minimax value $c$ such that $c< \frac dn S(a,b)^{\frac n{dp}}$ or $c< \frac dn S_R(a,b)^{\frac n{dp}}$, we obtain some existence results to (\ref{eq1.1}). We need some asymptotic estimates on the truncation function of the extremal function of $S_R(a,b)$. Let
$$
U_{\varepsilon}(x)=\frac1{(\varepsilon +|x|^{\frac{dp(n-p-pa)}{(p-1)(n-dp)}})^\frac{n-dp}{dp}},
$$

$$
k(\varepsilon)=c_0\big(\varepsilon(n-p-ap)\big)^\frac{n-dp}{dp}
$$
and $c_0$ is defined by (\ref{eq1.00013}). Then $y_\varepsilon(x):=k(\varepsilon) U_{\varepsilon}(x)$ is the extremal function of $S_R(a,b)$. Furthermore, there hold
\begin{equation}
\label{eq5.1}
\|Dy_\varepsilon\|^p_{L^p(\mathbb{R}^n, |x|^{-ap})}=S_R(a,b)^{\frac q{q-p}}=k(\varepsilon)^p \|DU_\varepsilon\|^p_{L^p(\mathbb{R}^n, |x|^{-ap})}\end{equation}
and
\begin{equation}
\label{eq5.2}
\|y_\varepsilon\|^q_{L^q(\mathbb{R}^n, |x|^{-bq})}=S_R(a,b)^{\frac q{q-p}}=k(\varepsilon)^q\|U_\varepsilon\|^q_{L^q(\mathbb{R}^n, |x|^{-bq})}.
\end{equation}
Let $\Omega\subset \mathbb{R}^n$ be an open bounded domain with $C^1$ boundary and $0\in \Omega$, $R>0$ such that $B_{2R}\subset \Omega$. Denote $u_{\varepsilon}(x)=\psi(x)U_{\varepsilon}(x)$ where $\psi(x)\equiv 1$ for $|x|<R$ and $\psi(x)\equiv 0$ for $|x|\geq 2R$. As $\varepsilon \to 0$, the behavior of $u_{\varepsilon}$ has to be the same as that of $U_{\varepsilon}$.

\begin{lemma}
\label{lem5.1}
Assume $1<p<n,\ -\infty< a<\frac{n-p}p,\ a\leq b\leq a+1,\ q=p^*(a,b)=\frac{np}{n-dp},\ d=1+a-b\in [0,\ 1],\ c>0$. Let
$$
v_\varepsilon(x)=\frac{u_{\varepsilon}(x)}{\|u_\varepsilon\|_{L^q(\Omega, |x|^{-bq})}}.
$$
Then $\|v_\varepsilon\|^q_{L^q(\Omega, |x|^{-bq})}=1$. Furthermore, there hold
\begin{enumerate}
\item[1.] $\|Dv_\varepsilon\|^p_{L^p(\Omega, |x|^{-ap})}= S_R(a,b) +O(\varepsilon^{(n-dp)/d})$;

\item[2.] $\|Dv_\varepsilon\|^\alpha_{L^\alpha(\Omega, |x|^{-ap})}= O(\varepsilon^{\frac{\alpha (n-dp)}{dp}})$ for $\alpha=1,2,p-2, p-1$;

\item[3.] $\|v_\varepsilon\|^p_{L^p(\Omega, |x|^{-(a+1)p+c})}=\begin{cases} O(\varepsilon^{(n-dp)/d})\ \mbox{ if } c>(n-p-ap)/(p-1)\\[2mm]
 O(\varepsilon^{(n-dp)/d}|\log \varepsilon|)\ \mbox{ if } c=(n-p-ap)/(p-1)\\[2mm]
 O(\varepsilon^{\frac{(p-1)(n-dp)(n+c-(a+1)p)}{dp(n-p-ap)} })\ \mbox{ if } c<(n-p-ap)/(p-1)

\end{cases}$
\end{enumerate}
\end{lemma}
The proof of Lemma \ref{lem5.1} is given in the Appendix B.

In the case where $a\geq 0, 1<p<n$, the results in \cite{HT} and \cite{CC} show that the minimizers of $S(a,b)$ are symmetric and given by (\ref{eq1.00013}). Combining Theorem \ref{thm4.2} and Lemma \ref{lem5.1}, there is the following existence result:
\begin{theorem}
\label{thm5.2} Let $\Omega\subset \mathbb{R}^n$ be an open bounded
domain with $C^1$ boundary and $0\in \Omega$, $1<p<n,\ 0\leq
a<\frac{n-p}p,\ a\leq b\leq a+1,\ q=p^*(a,b)=\frac{np}{n-dp},\
d=1+a-b\in (0,\ 1],\ c\leq (n-p-ap)/(p-1)$, and
$0<\lambda<\lambda_{1}$. Then there exists a nontrivial solution
$u\in {\cal D}_{a}^{1,p}(\Omega)$ to problem (\ref{eq1.1}).
\end{theorem}
\proof It is trivial that functional
$$
E_\lambda(u)=\frac1p\int_{\Omega} |x|^{-ap}|Du|^{p}\,dx-\frac1q \int_{\Omega} |x|^{-bq}|u|^{q}\,dx-\frac\lambda p \int_{\Omega} |x|^{-(a+1)p+c}|u|^p\,dx
$$
satisfies the geometric condition of the Mountain Pass Lemma
without (PS) condition due to  Ambrosetti and Rabinowitz
\cite{AR}. From Theorem \ref{thm4.2}, it suffices to show that
there exists a minimax value $c< \frac dn S(a,b)^{\frac n{dp}}$.
In fact, we will show that $\max\limits_{t\geq 0}
E_\lambda(tv_\varepsilon) <\frac dn S(a,b)^{\frac n{dp}}$ for
$\varepsilon$ small enough. Let
$$
\begin{array}{ll}
g(t)& =E_\lambda(tv_\varepsilon)\\[3mm]
& =\dfrac{t^p}p\displaystyle \int_{\Omega} |x|^{-ap}|Dv_\varepsilon|^{p}\,dx-\dfrac {t^q}q \displaystyle \int_{\Omega} |x|^{-bq}|v_\varepsilon|^{q}\,dx-\dfrac{\lambda t^p} p \displaystyle \int_{\Omega} |x|^{-(a+1)p+c}|v_\varepsilon|^p\,dx\\[3mm]
& =\dfrac{t^p}p\displaystyle \int_{\Omega} |x|^{-ap}|Dv_\varepsilon|^{p}\,dx-\dfrac {t^q}q -\dfrac{\lambda t^p} p \displaystyle \int_{\Omega} |x|^{-(a+1)p+c}|v_\varepsilon|^p\,dx.
\end{array}
$$
Since $0<\lambda<\lambda_{1}$, there holds $g(t)>0$ when $t$ is close to $0$, and $\lim\limits_{t\to \infty}g(t)=-\infty$ if $d=1+a-b\in (0,\ 1], q=p^*(a,b)=\frac{np}{n-dp}>p$. Thus $g(t)$ attains its maximum at some $t_\varepsilon>0$. From
$$
g^\prime(t)=t^{p-1}\big( \int_{\Omega} |x|^{-ap}|Dv_\varepsilon|^{p}\,dx -t^{q-p}-\lambda  \int_{\Omega} |x|^{-(a+1)p+c}|v_\varepsilon|^p\,dx\big)=0,
$$
there hold
$$
t_\varepsilon=\Big( \int_{\Omega} |x|^{-ap}|Dv_\varepsilon|^{p}\,dx-\lambda  \int_{\Omega} |x|^{-(a+1)p+c}|v_\varepsilon|^p \,dx\Big)^{1/(q-2)}
$$
and
$$
\begin{array}{ll}
g(t_\varepsilon)&=(\dfrac1p-\dfrac1q)\Big( \displaystyle\int_{\Omega} |x|^{-ap}|Dv_\varepsilon|^{p}\,dx-\lambda  \displaystyle\int_{\Omega} |x|^{-(a+1)p+c}|v_\varepsilon|^p \,dx\Big)^{q/(q-2)}\\[3mm]
& =\begin{cases}
\frac dn S(a,b)^{\frac n{dp}}+ O(\varepsilon^{\frac{n-dp}d})-O(\varepsilon^{\frac{(p-1)(n-dp)(n -(a+1)p+c)}{dp(n-p-ap)}} )\ \mbox{ if } c<\frac{n-p-ap}{p-1}\\[2mm]
 \frac dn S(a,b)^{\frac n{dp}}+ O(\varepsilon^{\frac{n-dp}d})-O(\varepsilon^{(n-dp)/d}|\log \varepsilon|)\ \mbox{ if } c=\frac{n-p-ap}{p-1}.
\end{cases}
\end{array}
$$
Note that for $c<(n-p-ap)/(p-1)$, there holds $\frac{n-dp}d > \frac{(p-1)(n-dp)(n -(a+1)p+c)}{dp(n-p-ap)}$. Thus for $\varepsilon$ small enough, there holds that $g(t_\varepsilon) <\frac dn S(a,b)^{\frac n{dp}}$.
\epf

In the case where $p=2$, combining Theorem \ref{thm4.3} and Lemma \ref{lem5.1}, there is the following existence result:

\begin{theorem}
\label{thm5.3}
Let $\Omega=B_1(0)$ is the unit ball in $\mathbb{R}^n$, $-\infty< a<\frac{n-2}2,\ a\leq b\leq a+1,\ q=2^*(a,b)=\frac{2n}{n-2d},\ d=1+a-b\in (0,\ 1],\ c\leq (n-2-2a)$, and $0<\lambda<\lambda_{1}$. Then there exists a nontrivial radial solution $u\in {\cal D}_{a,R}^{1,2}(\Omega)$ to problem (\ref{eq1.1}).
\end{theorem}
\proof It is trivial that functional
$$
E_\lambda(u)=\frac12\int_{\Omega} |x|^{-2a}|Du|^{2}\,dx-\frac1q \int_{\Omega} |x|^{-bq}|u|^{q}\,dx-\frac\lambda 2 \int_{\Omega} |x|^{-2(a+1)+c}|u|^2\,dx
$$
satisfies the geometric condition of the Mountain Pass Lemma
without (PS) condition due to  Ambrosetti and Rabinowitz
\cite{AR}. From Theorem \ref{thm4.3}, it suffices to show that
there exist a minimax value $c< \frac dn S_R(a,b)^{\frac n{2d}}$.
In fact, the same process in Theorem \ref{thm5.2} shows that $\max\limits_{t\geq 0}
E_\lambda(tv_\varepsilon) <\frac dn S_R(a,b)^{\frac n{2d}}$ for
$\varepsilon$ small enough for $c\leq n-2-2a$.
\epf

From the result in \cite{CC}, that is, $S(a,b)=S_R(a,b)$ for $p=2, a\geq 0$, Theorem \ref{thm4.2} and the proofs of Lemma \ref{lem5.1} and Theorem \ref{thm5.2} imply that
\begin{corollary}
\label{cor5.3}
Let $\Omega\subset \mathbb{R}^n$ be an open bounded domain with $C^1$ boundary and $0\in \Omega$, $0\leq a<\frac{n-2}2,\ a\leq b\leq a+1,\ q=2^*(a,b)=\frac{2n}{n-2d},\ d=1+a-b\in (0,\ 1],\ c\leq (n-2-2a)$, and $0<\lambda<\lambda_{1}$. Then there exists a nontrivial solution $u\in {\cal D}_{a}^{1,2}(\Omega)$ to problem (\ref{eq1.1}).
\end{corollary}

\begin{remark}
The results for the case where $a\geq 0, p=2$ had been obtained in \cite{CG} and \cite{NL} for $a=0, p=2$. But the results for the cases where $a<0$ or $p\neq 2$ had not been covered there.
\end{remark}

\begin{appendix}
\section{Appendix A}
\indent\textbf{Proof of Theorem \ref{thm3.2}. }Let $\{g_\varepsilon\}$ be a sequence of $C^2(\bar \Omega \setminus \{0\})$ functions converging to $g(\cdot, u)$ as $\varepsilon$ goes to $0^+$ and $u_\varepsilon$ the solution of
\begin{equation}
\label{eq3.6}
\left\{
\begin{array}{l}
-\mbox{div\,}(|x|^{-ap}(\varepsilon+|Du_\varepsilon|^{2})^{(p-2)/2}Du_\varepsilon)=g_\varepsilon,\mbox{ in } \Omega\\[2mm]
u_\varepsilon= 0, \ \ \mbox{on } \partial\Omega,
\end{array}
\right.
\end{equation}
Then from the standard regularity results in \cite{TP}, $u_\varepsilon$ is of class $C^3(\bar \Omega \setminus \{0\})$ and converges to $u$ in $C^{1, \alpha}(\bar \Omega \setminus \{0\})$, for some $\alpha\in (0,1)$. For problem (\ref{eq3.6}), we apply the Pohozaev integral identity-Lemma \ref{lem3.1} in $\Omega_\delta=\Omega\setminus \overline{B_\delta(0)}, 0<\delta <\mbox{dist\,}(0,\partial \Omega)$, noting that $u_\varepsilon$ may not vanish on the boundary $\partial B_\delta(0)=\{x\in \mathbb{R}^n\ :\ |x|=\delta\}$, or deduce directly by multiplying (\ref{eq3.6}) by $(A u_\varepsilon-h \cdot Du_\varepsilon)$ with $A=\frac np-(1+a)$, $h=x$, there holds
\begin{equation}
\label{eq3.7}
\begin{array}{ll}
&-\displaystyle \int_{\Omega_\delta} \mbox{div\,}(|x|^{-ap}(\varepsilon+|Du_\varepsilon|^{2})^{(p-2)/2}Du_\varepsilon)(A u_\varepsilon-x \cdot Du_\varepsilon)\,dx\\[3mm]
&\ \ \ \ \ \  =\displaystyle \int_{\Omega_\delta} g_\varepsilon(A u_\varepsilon-x \cdot Du_\varepsilon)\,dx.
\end{array}
\end{equation}
Integrating by parts over $\Omega_\delta$, there hold
\begin{equation}
\label{eq3.8}
\begin{array}{ll}
LHS&=-\displaystyle \int_{\partial \Omega_\delta}|x|^{-ap}(\varepsilon+|Du_\varepsilon|^{2})^{(p-2)/2}(A u_\varepsilon-x \cdot Du_\varepsilon)(Du_\varepsilon\cdot \nu) \,d\sigma\\[3mm]
&\ \ \ \ \ \ +  \displaystyle \int_{\Omega_\delta} |x|^{-ap}(\varepsilon+|Du_\varepsilon|^{2})^{(p-2)/2}Du_\varepsilon\cdot D( A u_\varepsilon-x \cdot Du_\varepsilon)\,dx\\[3mm]
&=-A \displaystyle \int_{|x|=\delta}|x|^{-ap}(\varepsilon+|Du_\varepsilon|^{2})^{(p-2)/2} u_\varepsilon (Du_\varepsilon\cdot \nu) \,d\sigma\\[3mm]
&\ \ \ \ \ \ +  \displaystyle \int_{\partial \Omega} |x|^{-ap}(\varepsilon+|Du_\varepsilon|^{2})^{(p-2)/2}|Du_\varepsilon|^{2}(x \cdot \nu)\,d\sigma\\[3mm]
&\ \ \ \ \ \ +  \displaystyle \int_{|x|=\delta} |x|^{-ap}(\varepsilon+|Du_\varepsilon|^{2})^{(p-2)/2}|Du_\varepsilon|^{2}(x \cdot \nu)\,d\sigma\\[3mm]
&\ \ \ \ \ \ + A \displaystyle \int_{\Omega_\delta} |x|^{-ap}(\varepsilon+|Du_\varepsilon|^{2})^{(p-2)/2}|Du_\varepsilon |^2\,dx\\[3mm]
&\ \ \ \ \ \ - \displaystyle \int_{\Omega_\delta} |x|^{-ap}(\varepsilon+|Du_\varepsilon|^{2})^{(p-2)/2}Du_\varepsilon\cdot D(x\cdot Du_\varepsilon)\,dx.
\end{array}
\end{equation}
Since $Du_\varepsilon\cdot D(x\cdot Du_\varepsilon)=|Du_\varepsilon|^{2}+\frac12 (x\cdot D(|Du_\varepsilon|^{2}))$, from (\ref{eq3.6}), there hold
\begin{equation}
\label{eq3.9}
\begin{array}{ll}
&\displaystyle \int_{\Omega_\delta} |x|^{-ap}(\varepsilon+|Du_\varepsilon|^{2})^{(p-2)/2}|Du_\varepsilon |^2\,dx\\[3mm]
&\ \ \ \ \ =\displaystyle \int_{\Omega_\delta} g_\varepsilon u_\varepsilon\,dx+\displaystyle \int_{|x|=\delta} |x|^{-ap}(\varepsilon+|Du_\varepsilon|^{2})^{(p-2)/2}u_\varepsilon (Du_\varepsilon\cdot \nu) \,d\sigma
\end{array}
\end{equation}
and
\begin{equation}
\label{eq3.10}
\begin{array}{ll}
&\dfrac12 \displaystyle \int_{\Omega_\delta} |x|^{-ap}(\varepsilon+|Du_\varepsilon|^{2})^{(p-2)/2} (x\cdot D(|Du_\varepsilon|^{2}))\,dx\\[3mm]
&\ \ \ \ \ \ \  =\dfrac1p \displaystyle \int_{\Omega_\delta}|x|^{-ap} x\cdot D((\varepsilon+|Du_\varepsilon|^{2})^{p/2} )\,dx\\[3mm]
&\ \ \  \ \ \ \  =\dfrac1p \displaystyle \int_{\partial \Omega}|x|^{-ap} (\varepsilon+|Du_\varepsilon|^{2})^{p/2} (x\cdot \nu)\,d\sigma\\[3mm]
&\ \ \ \ \ \ \ \ \ \ \ \ \ \  +\dfrac1p \displaystyle \int_{|x|=\delta}|x|^{-ap} (\varepsilon+|Du_\varepsilon|^{2})^{p/2} (x\cdot \nu)\,d\sigma\\[3mm]
&\ \ \ \ \ \ \ \ \ \ \ \ \ \ -\dfrac1p(n-ap)\displaystyle \int_{\Omega_\delta}|x|^{-ap} (\varepsilon+|Du_\varepsilon|^{2})^{p/2}\,dx,
\end{array}
\end{equation}
where $\nu$ is the unit outer normal vector. Substituting (\ref{eq3.9}) and (\ref{eq3.10}) into (\ref{eq3.8}) implies that
\begin{equation}
\label{eq3.11}
\begin{array}{ll}
LHS&=\displaystyle \int_{\partial \Omega} |x|^{-ap}(\varepsilon+|Du_\varepsilon|^{2})^{(p-2)/2}|Du_\varepsilon|^{2}(x \cdot \nu)\,d\sigma\\[3mm]
&\ \ \ \ \ \ +  \displaystyle \int_{|x|=\delta} |x|^{-ap}(\varepsilon+|Du_\varepsilon|^{2})^{(p-2)/2}|Du_\varepsilon|^{2}(x \cdot \nu)\,d\sigma\\[3mm]
&\ \ \ \ \ \ -\dfrac1p \displaystyle \int_{\partial \Omega}|x|^{-ap} (\varepsilon+|Du_\varepsilon|^{2})^{p/2} (x\cdot \nu)\,d\sigma\\[3mm]
&\ \ \ \ \ \  -\dfrac1p \displaystyle \int_{|x|=\delta}|x|^{-ap} (\varepsilon+|Du_\varepsilon|^{2})^{p/2} (x\cdot \nu)\,d\sigma\\[3mm]
&\ \ \ \ \ \ +(A-1) \displaystyle \int_{\Omega_\delta} g_\varepsilon u_\varepsilon\,dx\\[3mm]
&\ \ \ \ \ \
+\dfrac1p(n-ap)\displaystyle \int_{\Omega_\delta}|x|^{-ap} (\varepsilon+|Du_\varepsilon|^{2})^{p/2}\,dx.
\end{array}
\end{equation}
On the other hand, there holds
\begin{equation}
\label{eq3.12}
RHS=A \int_{\Omega_\delta} g_\varepsilon u_\varepsilon\,dx-\int_{\Omega_\delta} g_\varepsilon x \cdot Du_\varepsilon\,dx.
\end{equation}
Letting $\varepsilon\to 0^+$, there hold
\begin{equation}
\label{eq3.13}
\begin{array}{ll}
LHS&=(1-\dfrac1p)\displaystyle \int_{\partial \Omega} |x|^{-ap}|Du|^p (x \cdot \nu)\,d\sigma + (1-\dfrac1p) \displaystyle \int_{|x|=\delta} |x|^{-ap}|Du|^p (x \cdot \nu)\,d\sigma\\[3mm]
&\ \ \ \ \ \ +(A-1) \displaystyle \int_{\Omega_\delta} gu\,dx
+\dfrac1p(n-ap)\displaystyle \int_{\Omega_\delta}|x|^{-ap} |Du|^p \,dx
\end{array}
\end{equation}
and
\begin{equation}
\label{eq3.14}
\begin{array}{ll}
RHS&=A \displaystyle\int_{\Omega_\delta} g u\,dx-\displaystyle\int_{\Omega_\delta} gx \cdot Du\,dx\\[3mm]
&=A \displaystyle\int_{\Omega_\delta} g u\,dx-\displaystyle\int_{\partial \Omega_\delta} G(x,u)(x\cdot \nu)\,d\sigma \\[3mm]
&\ \ \ \ \ \ +\displaystyle \int_{\Omega_\delta}(x\cdot G_x)\,dx+n \displaystyle \int_{\Omega_\delta} G(x,u)\,dx.
\end{array}
\end{equation}
From (\ref{eq3.13}) and (\ref{eq3.14}), noting that $G(x,u)=\dfrac1q|x|^{-bq}|u|^q+\dfrac\lambda p |x|^{-p(1+a)+c}|u|^p$, there holds
\begin{equation}
\label{eq3.15}
\begin{array}{ll}
&(1-\dfrac1p)\displaystyle \int_{\partial \Omega} |x|^{-ap}|Du|^p (x \cdot \nu)\,d\sigma + (1-\dfrac1p) \displaystyle \int_{|x|=\delta} |x|^{-ap}|Du|^p (x \cdot \nu)\,d\sigma\\[3mm]
&\ \ \ \ \ \ \ \  +\dfrac1p(n-ap)\displaystyle \int_{\Omega_\delta}|x|^{-ap} |Du|^p \,dx \\[3mm]
&\ \ \  =\displaystyle\int_{\Omega_\delta} g u\,dx-\dfrac1q  \displaystyle \int_{|x|=\delta}|x|^{-bq}|u|^q (x \cdot \nu)\,d\sigma  -\dfrac\lambda p \displaystyle \int_{|x|=\delta}|x|^{-p(1+a)+c}|u|^p (x \cdot \nu)\,d\sigma\\[3mm]
&\ \ \ \ \ \ \ \  +(\dfrac nq-b) \displaystyle \int_{\partial \Omega}|x|^{-bq}|u|^q\,dx  + \lambda \dfrac{n-p(1+a)+c}{p} \displaystyle \int_{\partial \Omega}|x|^{-p(1+a)+c}|u|^p\,dx.
\end{array}
\end{equation}
Next, we need to get rid of the boundary integrals along $|x|=\delta$ in (\ref{eq3.15}). In fact, let $u$ be a solution of (\ref{eq1.1}), from the Caffarelli-Kohn-Nirenberg inequality (\ref{eq1.2}) or (\ref{eq1.4}), and Theorem \ref{thm2.1}, we know that
$$
\int_{\Omega} |x|^{-ap}|Du|^p\,dx,\quad \int_{\Omega} |x|^{-bq}|u|^q\,dx \mbox{  and  } \int_{\Omega} |x|^{-p(1+a)+c}|u|^p\,dx
$$
are finite. Therefore, by the mean-value theorem there exists a sequence $\{\delta_m\}, \delta_m\to 0^+$ such that integrals
$$
\int_{|x|=\delta} |x|^{-ap}|Du|^p (x \cdot \nu)\,d\sigma,
\int_{|x|=\delta}|x|^{-bq}|u|^q (x \cdot \nu)\,d\sigma, \int_{|x|=\delta}|x|^{-p(1+a)+c}|u|^p (x \cdot \nu)\,d\sigma \to 0
$$
as $m\to \infty$. Thus, letting $m\to \infty$ and noting (\ref{eq3.7}), we obtain (\ref{eq3.5}) from (\ref{eq3.15}).
\epf
\section{Appendix B}
\indent\textbf{Proof of Lemma \ref{lem5.1}. 1. }It is easy to see that
$$
\begin{array}{ll}
Du_\varepsilon(x)& =\begin{cases} DU_\varepsilon(x) \ \mbox{ if } |x|<R,\\[2mm]
U_\varepsilon(x)D\psi(x)+\psi(x)DU_\varepsilon(x) \ \mbox{ if } R\leq |x|<2R\\[2mm]
0 \ \mbox{ if } |x|\geq 2R
\end{cases} \\[4mm]
&=\begin{cases} -\dfrac{n-p-ap}{p-1}\ \dfrac{x}{(\varepsilon +|x|^{\frac{dp(n-p-pa)}{(p-1)(n-dp)}})^\frac{n}{dp} |x|^{2-\frac{dp(n-p-ap)}{(p-1)(n-dp)}}} \ \mbox{ if } |x|<R,\\[3mm]
U_\varepsilon(x)D\psi(x)+\psi(x)DU_\varepsilon(x) \ \mbox{ if } R\leq |x|<2R\\[2mm]
0 \ \mbox{ if } |x|\geq 2R,
\end{cases}
\end{array}
$$
$$
\begin{array}{ll}
\displaystyle \int_\Omega \dfrac{|Du_\varepsilon|^p}{|x|^{ap}}\,dx& =O(1)+ \displaystyle \int_{|x|<R} \dfrac{|DU_\varepsilon|^p}{|x|^{ap}}\,dx\\[3mm]
&=O(1)+ \displaystyle \int_{\mathbb{R}^n} \dfrac{|DU_\varepsilon|^p}{|x|^{ap}}\,dx\\[3mm]
&=O(1)+S_R(a,b)^{\frac q{q-p}}k(\varepsilon)^{-p}
\end{array}
$$
and
$$
\int_\Omega \dfrac{|u_\varepsilon|^q}{|x|^{bq}}\,dx=O(1)+S_R(a,b)^{\frac q{q-p}}k(\varepsilon)^{-q}.
$$
Thus, there holds
$$
\begin{array}{ll}
\|Dv_\varepsilon\|^p_{L^p(\Omega, |x|^{-ap})}& = \dfrac{\|Du_\varepsilon\|^p_{L^p(\Omega, |x|^{-ap})}}{\|u_\varepsilon\|^{p}_{L^q(\Omega, |x|^{-bq})} } \\[2mm]
&= \dfrac{O(1)+S_R(a,b)^{\frac q{q-p}}k(\varepsilon)^{-p}
}{O(1)+S_R(a,b)^{\frac p{q-p}}k(\varepsilon)^{-p}}\\[2mm]
&=S_R(a,b)+O(k(\varepsilon)^{p})=S_R(a,b)+O(\varepsilon^{(n-dp)/d}).
\end{array}
$$
\indent\textbf{2. } A direct computation shows that
$$
\begin{array}{ll}
&\displaystyle \int_\Omega \dfrac{|Du_\varepsilon|^\alpha}{|x|^{ap}}\,dx =O(1)+ \displaystyle \int_{|x|<R} \dfrac{|DU_\varepsilon|^\alpha}{|x|^{ap}}\,dx\\[3mm]
&\ \ \ \ \ =O(1)+ \displaystyle \int_{|x|<R} \big(\dfrac{n-p-ap}{p-1}\big)^\alpha \ \dfrac{|x|^{\alpha-ap}}{(\varepsilon +|x|^{\frac{dp(n-p-pa)}{(p-1)(n-dp)}})^\frac{\alpha n}{dp} |x|^{\alpha(2-\frac{dp(n-p-ap)}{(p-1)(n-dp)})}}\,dx\\[3mm]
&\ \ \ \ \ =O(1)+ \omega_n \displaystyle \int_0^R  \big(\dfrac{n-p-ap}{p-1}\big)^\alpha \ \dfrac{r^{\alpha-ap+n-1-\alpha(2-\frac{dp(n-p-ap)}{(p-1)(n-dp)})}}{(\varepsilon +r^{\frac{dp(n-p-pa)}{(p-1)(n-dp)}})^\frac{\alpha n}{dp} }\,dr\\[3mm]
&\ \ \ \ \ \leq O(1)+ \omega_n\big(\dfrac{n-p-ap}{p-1}\big)^\alpha \displaystyle \int_0^R r^{\alpha-ap+n-1-\alpha(2-\frac{dp(n-p-ap)}{(p-1)(n-dp)}) -\frac {\alpha (n-p-ap)}{(p-1)(n-dp)} }\,dr
\end{array}
$$
and the order of $r$ in the integrand is
$$
\begin{array}{ll}
& \alpha-ap+n-1-\alpha(2-\dfrac{dp(n-p-ap)}{(p-1)(n-dp)}) -\dfrac {\alpha (n-p-ap)}{(p-1)(n-dp)} \\[2mm]
&\ \ \ =\dfrac{np-n+\alpha-\alpha n-ap^2+ap+\alpha a p}{p-1}-1>-1
\end{array}
$$
for $\alpha=1,2, p-2, p-1$. Thus
$$
\int_\Omega \dfrac{|Du_\varepsilon|^\alpha}{|x|^{ap}}\,dx =O(1)
$$
and
$$
\begin{array}{ll}
\|Dv_\varepsilon\|^\alpha_{L^\alpha(\Omega, |x|^{-ap})} &=\dfrac{\|Du_\varepsilon\|^\alpha_{L^\alpha(\Omega, |x|^{-ap})}}{{\|u_\varepsilon\|^{\alpha }_{L^q(\Omega, |x|^{-bq})} } }\\[3mm]
& =\dfrac{O(1)}{O(1)+S_R(a,b)^{\frac \alpha {q-p}}k(\varepsilon)^{-\alpha}}\\[3mm]
&=O(k(\varepsilon)^{\alpha})=O(\varepsilon^{\frac{\alpha (n-dp)}{dp}}).
\end{array}
$$
\indent\textbf{3. } If $c=(n-p-ap)/(p-1)$, then there holds
$$
\begin{array}{ll}
\displaystyle \int_\Omega|x|^{-(a+1)p+c} |u_\varepsilon|^p \,dx &=O(1)+\displaystyle \int_{|x|<R}\dfrac1{(\varepsilon +|x|^{\frac{dp(n-p-pa)}{(p-1)(n-dp)}})^\frac{n-dp}{d} |x|^{ (a+1)p-c}}  \,dx\\[3mm]
&= O(1)+\omega_n \displaystyle \int_0^R \dfrac{r^{n-1-(a+1)p+c } }{(\varepsilon +r^{\frac{dp(n-p-pa)}{(p-1)(n-dp)}})^\frac{n-dp}{d} }\,dr\\[3mm]
&= O(1)+\omega_n  \displaystyle \int_0^{R\varepsilon^{-\frac{(p-1)(n-dp)}{dp(n-p-pa)} } } \dfrac{r^{n-1-(a+1)p+c } }{(1 +r^{\frac{dp(n-p-pa)}{(p-1)(n-dp)}})^\frac{n-dp}{d} }\,dr\\[3mm]
& \leq O(1)+\omega_n  \displaystyle \int_0^{R\varepsilon^{-\frac{(p-1)(n-dp)}{dp(n-p-pa)} } } \frac1r\,dr\\[3mm]
&=O(1)+O(|\log \varepsilon|).
\end{array}
$$
Then there holds
$$
\begin{array}{ll}
\|v_\varepsilon\|^p_{L^p(\Omega, |x|^{-(a+1)p+c})}& = \dfrac{\|u_\varepsilon\|^p_{L^p(\Omega, |x|^{-(a+1)p+c})}}{{\|u_\varepsilon\|^{p}_{L^q(\Omega, |x|^{-bq})} } } \\[3mm]
&=\dfrac{O(1)+O(|\log \varepsilon|)}{ {O(1)+S_R(a,b)^{\frac p{q-p}}k(\varepsilon)^{-p}} }\\[3mm]
&=O(k(\varepsilon)^{p}|\log \varepsilon|)=O(\varepsilon^{(n-dp)/d}|\log \varepsilon|).
\end{array}
$$

If $c>(n-p-ap)/(p-1)$, then there holds
$$
\begin{array}{ll}
\displaystyle \int_\Omega|x|^{-(a+1)p+c} |u_\varepsilon|^p \,dx &=O(1)+\displaystyle \int_{|x|<R}\dfrac1{(\varepsilon +|x|^{\frac{dp(n-p-pa)}{(p-1)(n-dp)}})^\frac{n-dp}{d} |x|^{ (a+1)p-c}}  \,dx\\[3mm]
&= O(1)+\omega_n \displaystyle \int_0^R \dfrac{r^{n-1-(a+1)p+c } }{(\varepsilon +r^{\frac{dp(n-p-pa)}{(p-1)(n-dp)}})^\frac{n-dp}{d} }\,dr\\[3mm]
& \leq O(1)+\omega_n  \displaystyle \int_0^{R}r^{n-1 -(a+1)p+c -\frac{p(n-p-ap)}{p-1}}\,dr\\[3mm]
&=O(1),
\end{array}
$$
the last equality is due to that $n-1 -(a+1)p+c -\frac{p(n-p-ap)}{p-1}>-1$ if $c>(n-p-ap)/(p-1)$. Thus there holds
$$
\begin{array}{ll}
\|v_\varepsilon\|^p_{L^p(\Omega, |x|^{-(a+1)p+c})}& = \dfrac{\|u_\varepsilon\|^p_{L^p(\Omega, |x|^{-(a+1)p+c})}}{{\|u_\varepsilon\|^{p}_{L^q(\Omega, |x|^{-bq})} } } \\[3mm]
&=\dfrac{O(1)}{ {O(1)+S_R(a,b)^{\frac p{q-p}}k(\varepsilon)^{-p}} }\\[3mm]
&=O(k(\varepsilon)^{p})=O(\varepsilon^{(n-dp)/d}).
\end{array}
$$

If $c< (n-p-ap)/(p-1)$, then $-\frac{n-dp}{d} +(n -(a+1)p+c) \frac{(p-1)(n-dp)}{dp(n-p-ap)}<0$ and $n-1 -(a+1)p+c -\frac{p(n-p-ap)}{p-1}<-1$, thus there hold
$$
\begin{array}{ll}
&\displaystyle \int_\Omega|x|^{-(a+1)p+c} |u_\varepsilon|^p \,dx=O(1)+\displaystyle \int_{|x|<R}\dfrac1{(\varepsilon +|x|^{\frac{dp(n-p-pa)}{(p-1)(n-dp)}})^\frac{n-dp}{d} |x|^{ (a+1)p-c}}  \,dx\\[4mm]
&\ \ \ \ = O(1)+\omega_n  \varepsilon^{-\frac{n-dp}{d} +(n -(a+1)p+c) \frac{(p-1)(n-dp)}{dp(n-p-ap)}} \displaystyle \int_1^\infty \dfrac{r^{n-1-(a+1)p+c } }{(1+r^{\frac{dp(n-p-pa)}{(p-1)(n-dp)}})^\frac{n-dp}{d} }\,dr\\[3mm]
&\ \ \ \ =O(\varepsilon^{-\frac{n-dp}{d} +(n -(a+1)p+c) \frac{(p-1)(n-dp)}{dp(n-p-ap)}} ),
\end{array}
$$
and
$$
\begin{array}{ll}
\|v_\varepsilon\|^p_{L^p(\Omega, |x|^{-(a+1)p+c})}& = \dfrac{\|u_\varepsilon\|^p_{L^p(\Omega, |x|^{-(a+1)p+c})}}{{\|u_\varepsilon\|^{p}_{L^q(\Omega, |x|^{-bq})} } } \\[3mm]
&=\dfrac{O(\varepsilon^{-\frac{n-dp}{d} +(n -(a+1)p+c) \frac{(p-1)(n-dp)}{dp(n-p-ap)}} )}{ {O(1)+S_R(a,b)^{\frac p{q-p}}k(\varepsilon)^{-p}} }\\[3mm]
&=O(\varepsilon^{\frac{(p-1)(n-dp)(n -(a+1)p+c)}{dp(n-p-ap)}} ).
\end{array}
$$
\epf
\end{appendix}
\end{document}